\documentclass[10pt]{amsart}
 \usepackage{latexsym,amsmath}
 \usepackage{amsmath,amsthm}
 \usepackage{amsfonts}
 \usepackage[psamsfonts]{amssymb}

\makeatletter
\@addtoreset{equation}{section}

\makeatother

 \newcommand{\beqa}{\begin{eqnarray}}
 \newcommand{\eeqa}{\end{eqnarray}}
 \newcommand{\bseq}{\begin{subequations}}
 \newcommand{\eseq}{\end{subequations}}

\font\tensmc=cmcsc10
\newcommand{\smc}{\tensmc}

\newcommand{\bbox}{\quad\hbox{\vrule \vbox{\hrule \vskip2pt \hbox{\hskip2pt
\vbox{\hsize=1pt}\hskip2pt} \vskip2pt\hrule}\vrule}}
\newcommand{\lessim}{\ \lower4pt\hbox{$
\buildrel{\displaystyle <}\over\sim$}\ }
\newcommand{\gessim}{\ \lower4pt\hbox{$\buildrel{\displaystyle >}
\over\sim$}\ }
\newcommand{\n}{\noindent}
\begin{document}
\title{\Large Exponential and Moment Inequalities for U-statistics}
\author{Evarist Gin\'e$^*$, Rafa{\l} Lata\l a$^\dagger$ and
Joel Zinn}

\maketitle
\pagestyle{myheadings}
\markboth{Inequalities for U-statistics}
{Gin\'e, Lata\l a and Zinn}

\renewcommand{\thefootnote}{\fnsymbol{footnote}}

\footnotetext[1]
{Research partially supported
 by NSF Grant No. DMS-96--25457.}

\footnotetext[2]
{Research partially supported by Polish Grant KBN 2
PO3A 043
15.}

\begin{abstract}
A Bernstein-type
exponential inequality for (generalized) canonical $U$-statistics of
order 2 is
obtained and the Rosenthal and Hoff-mann-J\o rgensen inequalities for sums of
independent random variables are extended to (generalized)
$U$-statistics of any
order whose kernels are either nonnegative or canonical.
\end{abstract}

\section{Introduction} Exponential inequalities, such as Bernstein's and
Prohorov's, and moment inequalities, such as Rosenthal's and
Hoffmann-J\o rgensen's,
are among the most basic tools for the analysis of sums of independent random
variables. Our object here consists in developing analogues of
such inequalities for
generalized $U$-statistics, in particular, for $U$-statistics and for
multilinear forms in independent random variables.

Hoffmann-J\o rgensen type moment inequalities for canonical (that is,
completely
degenerate) $U$-statistics of any order $m$ were first considered by
Gin\'e and Zinn (1992), and their version for $U$-statistics with
nonnegative kernels turned out to be useful for obtaining
best possible necessary integrability conditions in limit  theorems for
$U$-statistics. (By Khinchin's inequality it is irrelevant whether one
considers canonical or nonnegative kernels in moment inequalities, at least if
multiplicative  constants are not at
issue). Klass and Nowicki (1997) also obtained moment inequalities for
nonnegative generalized
$U$-statistics, but only for order $m=2$, and their decomposition of
the moments is
more complete than that in Gin\'e and Zinn (1992). Ibragimov and Sharakhmetov
(1998, 1999) recently obtained analogues of Rosenthal's inequality for
nonnegative and for canonical
$U$-statistics. The moment inequalities we present in the first part of this
article,  valid for canonical and for nonnegative generalized
$U$-statistics of any order $m$, when  specialized
to $m=2$, represent the same level of moment decomposition as the
Klass-Nowicki
inequalities, coincide with theirs for powers $p>1$ (except for constants) and
are  expressed in terms of
different, simpler quantities for powers $p<1$. Proposition 2.1 below, which
constitutes the first step towards more elaborate bounds such as those in
Theorem 2.3 below, has also been obtained, up to constants, by Ibragimov and
Sharakhmetov.  Our  proofs consist of simple iterations of the classical
moment inequalities for sums of independent random variables.

The moment inequalities in the first part of this article do imply exponential
bounds for canonical $U$-statistics of any order and with bounded kernels which
are sharper than those in Arcones and Gin\'e (1993); however, they are not of
the best kind as they do not exhibit Gaussian behavior for part of the tail,
which they should in view of the tail behavior of Gaussian
chaos.

In the second part of this article we improve the moment inequalities from the
first part in the case of generalized canonical $U$-statistics of order 2, and
for moments of order $p\ge2$ (Theorem 3.2). The bounds not only involve moments
but also the
$L_2$ operator norm of the matrix of kernels. Then we show how these improved
moment inequalities imply what we believe is the correct analogue (up to
constants) of Bernstein's exponential inequality for generalized canonical
$U$-statistics of order 2 (Theorem 3.3). This exponential inequality, which
does exhibit  Gaussian behavior for small values of $t$, is strong enough to
imply the law of the iterated logarithm for canonical $U$-statistics under
conditions which are also necessary. The main new ingredient in this part of
the paper is Talagrand's (1996) exponential bound for empirical
processes, which gives a Rosenthal-Pinelis type inequality
for moments of empirical processes (Proposition 3.1) basic for the derivation
of the moment inequality for $U$-statistics of order 2.

Because of the decoupling results of de la Pe\~na and
Montgomery-Smith (1995), we
can work with decoupled $U$-statistics, and this allows us to
proceed by conditioning
and iteration.

\section{Moment inequalities} We consider estimation of moments of {\it
generalized decoupled $U$-statistics}, defined as
\begin{equation}
\sum_{1\le i_1,\dots,i_m\le
n}h_{i_1,\dots,i_m}(X_{i_1}^{(1)},\dots,X_{i_m}^{(m)}),
\end{equation}
where the random variables $X_i^{(j)}:1\le i\le n, 1\le j\le m$, $m\le n$, are
independent (not necessarily with the same distribution) and  take
values in a measurable space $(S,{\mathcal S})$, and $h_{i_1,\dots,i_m}$
are  real valued measurable functions on $S^m$.
For short, this sum is denoted by $\sum_{\bf i}h_{\bf i}$.

Given $J\subseteq\{1,\dots, m\}$ ($J=\emptyset$ is not excluded), and ${\bf
i}=(i_1,\dots, i_m)\in\{1,\dots,n\}^m$ we set ${\bf i}_J$ to be the point of
$\{1,\dots,n\}^{|J|}$ obtained from $\bf i$ by deleting the coordinates in the
places not in $J$ (e.g., if ${\bf i}=(3,4,2,1)$ then ${\bf
i}_{\{1,3\}}=(3,2)$).
Also,
$\sum_{{\bf i}_J}$ indicates sum over $1\le i_j\le n$, $j\in J$ (for
instance, if
$m=4$ and
$J=\{1,3\}$, then
$$\sum_{{\bf
i}_J}h_{\bf i}=\sum_{{\bf i}_{\{1,3\}}}h_{i_1,i_2,i_3,i_4}=\sum_{1\le i_1,i_3
\le
n}h_{i_1,i_2,i_3,i_4}(X_{i_1}^{(1)},\dots,X_{i_4}^{(4)}).)$$ By convention,
$\sum_{{\bf i}_\emptyset}a=a$.

Likewise, while $E$ will denote expected value with respect to all the
variables,
$E_J$ will denote expected value only with respect to the variables
$X_i^{(j)}$ with $j\in J$ and $i\in\{1,\dots,n\}$. By convention, $E_\emptyset
a=a$.

Rosenthal's inequality is easiest to extend to $U$-statistics because it
involves
only moments of sums (as opposed to moments of maxima and quantiles for
Hoffmann-J\o rgensen's inequality). So, we will first obtain analogues of
Rosenthal's inequality, and then we will transform these inequalities into
analogues of Hoffmann-J\o rgensen's by first showing that some
moments of sums can
be replaced by moments of maxima, and then, that the lowest
moment can in fact be
replaced by a quantile. We will illustrate this three-steps
procedure first in the
case of nonnegative kernels and moments of order $p\ge 1$.
Then we will see that
this also solves, via Khinchin's inequality, the case of canonical kernels and
moments of order $p\ge 2$. Finally, we will consider the
case of moments of order
$p<1$ for positive kernels and $p<2$ for canonical, cases
in which the inequalities
are less neat, but still useful. We will pay some attention to
the behavior of the
constants as $p\to\infty$ in these inequalities since such behavior
translates into
(exponential) integrability properties.

\vskip.1truein
\n{\it 2.1. Nonnegative kernels, moments of order $p\ge 1$.} For nonnegative
independent random
variables $\xi_i$, we have the following two improvements of Rosenthal's
inequalities, valid for
$p\ge 1$:

\n 1) Lata\l a's, 1997:
$$E\Bigl(\sum \xi_i\Bigr)^p\le (2e)^p\max\biggl[{\frac{e}{p}}p^p\sum
E\xi_i^p,~e^p\Bigl(\sum E\xi_i\Bigr)^p\biggr],~~p>1,\leqno(R_1)$$
(see Pinelis (1994) for the corresponding inequality when the random variables
are centered);

\n 2) Johnson, Schechtman and Zinn's, 1985:
$$E\Bigl(\sum \xi_i\Bigr)^p\le K^p\biggl(\frac{p}{\log p}\biggr)^p
\max\biggl[\sum
E\xi_i^p,~\Bigl(\sum E\xi_i\Bigr)^p\biggr],~~p>1,\leqno(R_2)$$
where $K$ is a universal constant. See Utev (1985) and Figiel, Hitczenko,
Johnson, Schechtman and Zinn (1997) for more precise inequalities of the same
type.

\n And for general $p>0$, we have the following improved Hoffmann-J\o rgensen
inequality, that follows from Kwapie\'n and Woyczy\'nski (1992) and which
can be
obtained as in the proof of Theorem 1.2.3 in de la Pe\~na and Gin\'e (1999):

\n 3)$$E\Big\|\sum\xi_i\Big\|^p\le 2^{p-2}\cdot
2^{(p-1)\vee0}\cdot(p+1)^{p+1}\Bigl[t_0^p+E\max\|\xi_i\|^p\Big],~~~p>0,
\leqno(H)$$
where
$$t_0:=\inf\biggl[t>0:\Pr\Bigl\{\Big\|\sum\xi_i\Big\|>t\Bigr\}\le \frac{1}
{2}\biggr],$$
and where we write norm for absolute value in order to include not only
independent
nonnnegative real random variables, but also independent nonnegative random
functions
$\xi_i$ taking values in certain `rearrangement invariant spaces' such as
$L_s(\Omega,
\Sigma,
\mu)$,  $0< s<\infty$, with
$\|\xi\|:=\bigl(\int|\xi|^sd\mu\bigr)^{1/(s\vee 1)}$, or $\ell_\infty(L_s)$.
Note that, by Markov,
$$t_0\le 2^{1/r}\biggl(E\Big\|\sum\xi_i\Big\|^r\biggr)^{1/r},$$
so that, $(H)$ becomes:

\n 4) for $0<r<p<\infty$,
$$E\Big\|\sum\xi_i\Big\|^p\le 2^{p-2}\cdot
2^{(p-1)\vee0}\cdot(p+1)^{p+1}\biggl[2^{p/r}
\biggl(E\big\|\sum\xi_i\big\|^r\biggr)^{p/r}$$
\vskip-.2truein
$$\mbox{\phantom{aaaaaaaaaaaaaaaaaaaaaaaaaaaaaaaaaaaaaa}}+E\max\|\xi_i\|^p\bigg]
\leqno(H_r)$$
Inequalities $(H)$ and $(H_r)$ hold for spaces of functions which are
quasinormed measurable linear spaces whose quasinorm $\|\cdot\|$ has the
property
 that $\|x\|\le \|y\|$ whenever
$0\le x\le y$.

 In
the following proposition we extend inequalities
$(R_1)$ and
$(R_2)$ by  means of
an easy induction.

\vskip.1truein
\n{\bf Proposition 2.1.} {\it Let $m\in\bf N$, $p>1$, and, for all
${\bf i}\in\{1,\dots,n\}^m$, let $h_{\bf
i}$ be a nonnegative function of $m$ variables whose $p$-th power is integrable
for the law of ${\bf X}_{\bf i}=(X_{i_1}^{(1)},\dots,X_{i_m}^{(m)})$. Then,}
\begin{eqnarray}
\lefteqn{\max_{J\subseteq\{1,\dots,m\}}\biggl[\sum_{{\bf i}_J}
E_J\Bigl(\sum_{{\bf
i}_{J^c}}E_{J^c}h_{\bf i}\Bigr)^p\biggr]
\le  E\Bigl(\sum_{\bf i}h_{\bf i}\Bigr)^p} \nonumber \\
&&\le
(2e^2)^{mp}\sum_{J\subseteq\{1,\dots,m\}}\biggl[p^{|J|p}\sum_{{\bf
i}_J}E_J\Bigl(\sum_{{\bf i}_{J^c}}E_{J^c}h_{\bf
i}\Bigr)^p\biggr],
\end{eqnarray}
{\it and also, there exists a
universal constant $K<\infty$ such that}
$$E\Bigl(\sum_{\bf i}h_{\bf i}\Bigr)^p
\le
K^{mp}\biggl(\frac{p}{\log
p}\biggr)^{mp}\max_{J\subseteq\{1,\dots,m\}}\biggl[\sum_{{\bf
i}_J}E_J\Bigl(\sum_{{\bf i}_{J^c}}E_{J^c}h_{\bf
i}\Bigr)^p\biggr].\leqno(2.2')$$

\vskip.05truein
\n {\bf Proof.} The proof of (2.2') with sum over the subsets $J$ instead of
maximum differs from that of (2.2) only in the starting point ($(R_2)$ instead
of $(R_1)$); then, replacing sum by maximum simply increases the constant by
a factor of $2^m$. The left side inequality in (2.2) follows by H\"older since
$p\ge 1$. Consider the right hand side inequality. For
$m=1$ this is just inequality ($R_1$) and we can proceed by induction.
Suppose the
result holds for
$m-1$.  By applying the induction hypothesis to
$$E\Bigl(\sum_{\bf i}h_{\bf i}\Bigr)^p=E_{m}E_{\{1,\dots,m-1\}}
\biggl[\sum_{{\bf
i}_{\{1,\dots,m-1\}}}\Big(\sum_{i_m}h_{\bf i}\Big)\biggr]^p,$$
 we only have to consider the
generic term in the decomposition (2.2) for the new kernels
$\Big(\sum_{i_m}h_{\bf i}\Big)$ with the $X_i^{(m)}$ variables fixed. In other
words,
letting
$J_{m-1}$ be any subset of $\{1,\dots,m-1\}$ and  $J_{m-1}^c$ its
complement with respect to $\{1,\dots,m-1\}$, we must estimate
\begin{eqnarray}\lefteqn{E_{m}\sum_{{\bf
i}_{J_{m-1}}}E_{J_{m-1}}\biggl(\sum_{{\bf
i}_{J_{m-1}^c}}E_{J_{m-1}^c}\Big(\sum_{i_m}h_{\bf i}\Big)\biggr)^p} \nonumber
\\ &&=\sum_{{\bf
i}_{J_{m-1}}}E_{J_{m-1}}E_{m}\biggl(\sum_{i_m}
\Big(E_{J_{m-1}^c}\sum_{{\bf i}_{J_{m-1}^c}}h_{\bf
i}\Big)\biggr)^p. \nonumber
\end{eqnarray}
Rosenthal's inequality ($R_1$) applied to the kernels
$E_{J_{m-1}^c}\sum_{{\bf i}_{J_{m-1}^c}}h_{\bf
i}$ with the variables in $J_{m-1}$ fixed, gives
\begin{eqnarray}E_{m}\biggl(\sum_{i_m}
\Big(E_{J_{m-1}^c}\sum_{i_{J_{m-1}^c}}h_{\bf i}\Big)\biggr)^p\!\!\!&\le&\!\!\!
(2e^2)^p\biggl[\Bigl(\sum_{i_m,i_{J_{m-1}^c}}E_{m}E_{J_{m-1}^c}
h_{\bf i}\Bigr)^p \nonumber
\\
&&\mbox{\phantom{aa}}+p^p\sum_{i_m}E_m\Bigl(E_{J_{m-1}^c}\sum_{{\bf
i}_{J_{m-1}^c}}h_{\bf i}\Bigr)^p\biggr].\nonumber
\end{eqnarray}
Upon integrating  each term with respect to $E_{J_{m-1}}$ and summing over
${\bf
i}_{J_{m-1}}$, we then obtain
\begin{eqnarray}\lefteqn{E_{m}\sum_{{\bf
i}_{J_{m-1}}}E_{J_{m-1}}\biggl(\sum_{{\bf
i}_{J_{m-1}^c}}E_{J_{m-1}^c}\Big(\sum_{i_m}h_{\bf i}\Big)\biggr)^p}\nonumber\\
&&\le
(2e^2)^p\biggl[\sum_{{\bf
i}_{J_{m-1}}}E_{J_{m-1}}\Bigl(\sum_{{\bf i}_{J_{m-1}^c\cup\{m\}}}
E_{J_{m-1}^c\cup\{m\}}h_{\bf
i}\Bigr)^p\nonumber \\
&&
\mbox{\phantom{aaa}}+p^p\sum_{{\bf i}_{J_{m-1}\cup\{m\}}}
E_{J_{m-1}\cup\{m\}}\Bigl(E_{J_{m-1}^c}\sum_{{\bf
i}_{J_{m-1}^c}}h_{\bf i}\Bigr)^p\biggr].\nonumber
\end{eqnarray}
Multiplying by $(2e^2)^{(m-1)p}p^{|J_{m-1}|}$, this is the sum of two terms of
the form
\hfil\break
$(2e)^{mp}p^{|J|p}\sum_{{\bf
i}_J}E_J\Bigl(\sum_{{\bf i}_{J^c}}E_{J^c}h_{\bf i}\bigr)^p$ (for $J=J_{m-1}$
and for
$J=J_{m-1}\cup\{m\}$), proving the proposition.\quad$\bbox$

\vskip.2truein
This proposition solves the problem of estimating, up to constants,
the moments of
a decoupled
$U$-statistic by `computable' expressions. For instance, if the
functions $h_{\bf
i}$ are all equal and if the variables $X_i^{(j)}$ are i.i.d., then the typical
term at the right of (2.1) just becomes $n^{|J|+p|J^c|}E_J(E_{J^c}h)^p$,
a `mixed
moment' of $h$. For $m=2$ the right hand side of inequality (2.2) is just:
\begin{eqnarray}
E\biggl(\sum_{i,j}h_{i,j}(X_i^{(1)},X_j^{(2)})\biggr)^p
\!\!\!&\le&\!\!\!(2e^2)^{2p}
\biggl[\Bigl(\sum_{i,j}Eh_{i,j}(X_i^{(1)},X_j^{(2)})\Bigr)^p \nonumber \\
&&\mbox{}+p^p\sum_iE_1\Bigl(\sum_jE_2h_{i,j}(X_i^{(1)},X_j^{(2)})\Bigr)^p
\nonumber \\
&&\mbox{\phantom{}}+p^p\sum_jE_2\Bigl(\sum_iE_1h_{i,j}(X_i^{(1)},X_j^{(2)})
\Bigr)^p \nonumber
\end{eqnarray}
\vskip-.2truein
$$\leqno(2.2'')\mbox{\phantom{aaaaaaaaaaaaaaaaaaaaaaaaa}}+p^{2p}
\sum_{i,j}Eh_{i,j}^p
(X_i^{(1)},X_j^{(2)})\biggr].
$$

We have been careful with the dependence on $p$ of the constants
because it is of some interest to obtain constants of the best order as
$p\to\infty$. In fact, (2.2') exhibits constants of the best order as
can be seen by
taking the product of two independent copies of the example in Johnson,
Schechtman and Zinn (1985), Proposition 2.9.

Next we replace the external sums of expected values at the right
side of the above
inequalities by expectations of maxima without significantly altering the order
of the multiplicative constants. If $\xi_i$ are independent nonnegative random
variables, then,
\begin{equation}
\frac{1}{2}\biggl[\delta_0^p\vee\sum
E\xi_i^pI_{\xi_i>\delta_0}\biggr]\le E\max\xi_i^p\le
\delta_0^p+\sum E\xi_i^pI_{\xi_i>\delta_0},~~0<p<\infty,
\end{equation}
where
\begin{equation}
\delta_0=\inf\biggl[t>0:\sum\Pr\bigl\{\xi_i>t\bigr\}\le
1\biggr]
\end{equation}
(Gin\'e and Zinn (1983); see also de la Pe\~na and Gin\'e (1999), page 22). The
left hand side of (2.3) gives  that, for
$0<r< p$ and $\xi_i$ independent,
\begin{equation}
\sum E|\xi_i|^p\le 2E\max|\xi_i|^p+2\Bigl(\sum E
|\xi_i|^r\Bigr)\Bigl(E\max|\xi_i|^p\Bigr)^{(p-r)/p}
\end{equation}
(e.g., de la Pe\~na and Gin\'e (1999), page
48).
This inequality, applied with $r=1<p$, yields
\begin{equation}
p^{\alpha p}\sum
E|\xi_i|^p\le2(1+p^\alpha)\max\biggl[p^{\alpha p}E\max|\xi_i|^p,~\Bigl(\sum
E|\xi_i|\Bigr)^p\biggr]
\end{equation}
for all  $\alpha\ge0$.
There are similar inequalities for other values
of $r$; $r=1$ is adequate for
$\xi_i\ge 0$, but $r=2$ is better for centered variables. If we use inequality
(2.6)  in (2.2''), iteratively for the last term, we obtain that, for a
universal
constant
$K$ (easy but cumbersome to compute),
$h_{i,j}\ge0$,
$p> 1$,
\begin{eqnarray}E\Bigl(\sum_{i,j}h_{i,j}\Bigr)^p\!\!\!&\le&\!\!\!
K^p(2e^2)^pp^4\biggl[\Bigl(\sum_{i,j}Eh_{i,j}\Bigr)^p
+p^pE_1\max_i\Bigl(\sum_jE_2h_{i,j}\Bigr)^p \nonumber \\
&&\mbox{\phantom{aaaaaaaaa}}+p^pE_2\max_j\Bigl(\sum_iE_1h_{i,j}\Bigr)^p+
p^{2p}E\max_{i,j}h_{i,j}^p\biggr].
\end{eqnarray}
Inequality (2.7) was obtained, up to constants, by Klass and Nowicki (1997) (it
is their inequality (4.14)). Our proof is different, and it is contained in the
proof of the next corollary, which extends inequality (2.7) to any
$m$.

\vskip.1truein
\n{\bf Corollary 2.2.} {\it Under the same hypotheses as in Proposition
2.1, there
exist universal constants $K_m$
such that}
\begin{eqnarray}\lefteqn{\max_{J\subseteq\{1,\dots,m\}}\biggl[
E_J\max_{{\bf i}_J}\Bigl(\sum_{{\bf
i}_{J^c}}E_{J^c}h_{\bf i}\Bigr)^p\biggr]
\le E\Bigl(\sum_{\bf i}h_{\bf i}\Bigr)^p} \nonumber \\
&&\le
K_m^p\sum_{J\subseteq\{1,\dots,m\}}\biggl[p^{|J|p}E_J\max_{{\bf
i}_J}\Bigl(\sum_{{\bf i}_{J^c}}E_{J^c}h_{\bf i}\Bigr)^p\biggr],
\end{eqnarray}
{\it and}
$$E\Bigl(\sum_{\bf i}h_{\bf i}\Bigr)^p
\le
K_m^p\biggl(\frac{p}{\log
p}\biggr)^{mp}\max_{J\subseteq\{1,\dots,m\}}\biggl[E_J\max_{{\bf
i}_J}\Bigl(\sum_{{\bf i}_{J^c}}E_{J^c}h_{\bf i}\Bigr)^p\biggr].\leqno(2.8')$$

\n {\bf Proof.} The left side of (2.8) follows by H\"older. Inequality (2.8')
has a proof similar to that of the right hand side of (2.8), and therefore we
only prove the latter. We will prove it by induction over $m$ simultaneously
with the inequality
\begin{equation}
\label{addin}
p^{mp}\sum_{{\bf i}}Eh_{{\bf i}}^{p}\leq \tilde{K}_{m}^{p}
\sum_{J\subset\{1,\ldots,m\}}
\biggl[p^{|J|p}E_{J}\max_{{\bf i}_{J}}\Bigl(\sum_{{\bf i}_{J^{c}}}E_{J^{c}}
h_{{\bf i}}\Bigr)^{p}\biggr].
\end{equation}
Let us first note that the inequalities (\ref{addin}) for $1,\ldots,m-1$
together with (2.2) imply (2.8). It is therefore enough to show that
if (2.8) and (\ref{addin}) hold for $1,\ldots,m-1$ then (\ref{addin})
is satisfied for $m$. We will follow the notation of the proof of Proposition
2.1. Inequality (2.9) for $m=1$ is just (2.6), and (2.8) for $m=1$ is
$(H_1)$ (which also follows from $(R_1)$ and (2.6)). By the induction
assumptions we have
\begin{equation}
\label{ind1}
p^{mp}\sum_{{\bf i}}Eh_{{\bf i}}^{p}=p^{p}\sum_{i_{m}}E_{m}
  p^{(m-1)p}\sum_{{\bf i}_{\{1,\ldots,m-1\}}}E_{\{1,\ldots,m-1\}}
  h_{{\bf i}}^{p}
\end{equation}
\vskip-.1truein
\n $\le\tilde{K}_{m-1}^{p}$
\vskip-.1truein
$$\times\sum_{J_{m-1}\subset\{1,\ldots,m-1\}}
\biggl[p^{(|J_{m-1}|+1)p}E_{J_{m-1}}E_{m}\sum_{i_{m}}\max_{{\bf i}_{J_{m-1}}}
\Bigl(\sum_{{\bf i}_{J_{m-1}^{c}}}E_{J_{m-1}^{c}}
h_{{\bf i}}\Bigr)^{p}\biggr].$$
Now, by (2.6), for any $J_{m-1}\subset\{1,\ldots,m-1\}$ we have
\begin{eqnarray}\label{ind2}\lefteqn{
p^{(|J_{m-1}|+1)p}E_{J_{m-1}}E_{m}\sum_{i_{m}}\max_{{\bf
i}_{J_{m-1}}} \Bigl(\sum_{{\bf i}_{J_{m-1}^{c}}}E_{J_{m-1}^{c}}
h_{{\bf i}}\Bigr)^{p}}\nonumber\\
&&\le2(1+p)\biggl[p^{(|J_{m-1}|+1)p}E_{J_{m-1}\cup \{m\}}
  \max_{{\bf i}_{J_{m-1}\cup\{m\}}}\Bigl(\sum_{{\bf
i}_{J_{m-1}^{c}}}E_{J_{m-1}^{c}}
  h_{{\bf i}}\Bigr)^{p}\nonumber\\
&&\mbox{\phantom{aaaaaaaaaa}}+p^{|J_{m-1}|p}E_{J_{m-1}}\Bigl(\sum_{i_{m}}
E_{m}\max_{{\bf
i}_{J_{m-1}}}
\sum_{{\bf i}_{J_{m-1}^{c}}}E_{J_{m-1}^{c}}h_{{\bf i}}\Bigr)^{p}\biggr].
\end{eqnarray}
To estimate the last term we note that
\begin{eqnarray}\label{ind3}
\lefteqn
{p^{|J_{m-1}|p}E_{J_{m-1}}\Bigl(\sum_{i_{m}}E_{m}\max_{{\bf i}_{J_{m-1}}}
\sum_{{\bf i}_{J_{m-1}^{c}}}E_{J_{m-1}^{c}}h_{{\bf i}}\Bigr)^{p}}\nonumber\\
&&\leq
p^{|J_{m-1}|p}E_{J_{m-1}}\Bigl(\sum_{{\bf i}}E_{J_{m-1}^{c}\cup\{m\}}
  h_{{\bf i}}\Bigr)^{p}\mbox{\phantom{aa}}\end{eqnarray}
\vskip-.1truein
$$\leq
\tilde{K}_{|J_{m-1}|}^{p}\sum_{J\subset J_{m-1}}\!p^{|J|p}E_{J}\max_{{\bf
i}_{J}}
\Bigl(\sum_{{\bf i}_{(J_{m-1}\setminus J)\cup J_{m-1}^{c}\cup\{m\}}}
\!\!\! E_{(J_{m-1}\setminus J)\cup J_{m-1}^{c}\cup\{m\}}h_{{\bf
i}}\Bigr)^{p},$$ where in the last line we use the induction assumption (2.8)
for $|J_{m-1}|<m$. Finally (\ref{ind1}), (\ref{ind2}) and (\ref{ind3}) imply
(\ref{addin}) and complete the proof. \quad$\bbox$

\vskip.1truein
\n{\bf Remark.} The proof of Proposition 2.6 below will use a version of
Corollary 2.2 for nonnegative random functions taking values in $L_r$. The
inequality is as follows: for $p>1$ there exists $K_{m,p,r}<\infty$ such
that
$$E\bigl\|\sum_{\bf i}h_{\bf i}\bigr\|^p
\le
K_{m,p,r}\max_{J\subseteq\{1,\dots,m\}}\biggl[E_J\max_{{\bf
i}_J}\Bigl(E_{J^c}\big\|\sum_{{\bf i}_{J^c}}h_{\bf
i}\big\|\Bigr)^p\biggr].\leqno(2.8'')$$
The proof is similar to the previous ones and is omited: one takes $(H_p)$ as
the starting point of the induction.

\vskip.1truein
Finally we come to the third step, which will extend
Hoffmann-J\o rgen-sen's inequality
$(H)$ for $p\ge 1$. If we want to use the inequalities from
Corollary 2.2 to obtain
boundedness of moments from stochastic boundedness of a sequence of
$U$-statistics, we need to replace the term corresponding to
$J=\emptyset$ by the
$p$-th power of a quantile of $\sum_{\bf i} h_{\bf i}$. For this we use
Paley-Zygmund's inequality (e.g., Kahane (1968) or de la Pe\~na
and Gin\'e (1999)):
if $A$ is a nonnegative
random variable and $0<r<p<\infty$, then, for all $0<\lambda<1$,
\begin{equation}
\Pr\Bigl\{A>\lambda\|A\|_r\Bigr\}\ge\biggl[(1-\lambda^r)
\frac{\|A\|_r}{\|A\|_p}\biggr]^{p/(p-r)},
\end{equation}
where $\|A\|_r=\bigl(E |A|^r\bigr)^{1/r}$ for $0<r<\infty$.
Consider for instance
inequality (2.8). It has the form
$$EA^p\le B+K_m^p(EA)^p, \ \ p>1,$$
with $A=\sum_{\bf i} h_{\bf i}$. Then, either $B\ge K_m^p(EA)^p$,
in which case
we have $EA^p\le 2B$, or $B< K_m^p(EA)^p$, in which case we have $EA^p\le
2K_m^p(EA)^p$ and we can apply Paley-Zygmund's (2.13) with
$\lambda=1/2$ and $r=1$. It gives
$$
\Pr\Bigl\{A>\frac{1}{2}EA\Bigr\}\ge\frac{1}{
2^{(p+1)/(p-1)}K_m^{p/(p-1)}}.$$ Hence, if we define
\begin{equation}
t_0=\inf\Bigl[t\ge 0: \Pr\{A>t\}\le\frac{1}{
2^{(p+1)/(p-1)}K_m^{p/(p-1)}}\Bigr],
\end{equation}
we obtain $EA\le 2t_0$. So, in either case,
$$EA^p\le 2B+2^{1+p}K_m^{p}t_0^p.$$
Also, by Markov's inequality,
$$\frac{1}{
2^{(p+1)/(p-1)}K_m^{p/(p-1)}}t_0^p\le EA^p.$$
We then have:

\vskip.1truein
\n{\bf Theorem 2.3.} {\it Under the hypotheses of Proposition 2.1,
there exist a
universal constants
$K_m<\infty$ such that, if $t_0$ is as defined by (2.14) for $A=\sum_{\bf i}
h_{\bf i}$, then}
\begin{eqnarray}\lefteqn{\frac{1}{(4K_m)^{p/(p-1)}}t_0^p\vee
\max\biggl[\max_{\begin{array}{c}\scriptstyle J\subseteq\{1,\dots,m\} \\
\scriptstyle J\ne\emptyset \end{array} }\biggl[E_J\max_{{\bf
i}_J} \Bigl(\sum_{{\bf
i}_{J^c}}E_{J^c}h_{\bf i}\Bigr)^p\biggr]}\nonumber \\
&&\le E\Bigl(\sum_{\bf i}h_{\bf i}\Bigr)^p\nonumber \\
&&\le
(4K_m)^{p}\Biggl\{2^{1+p}t_0^p+\sum_{\begin{array}{c}\scriptstyle
J\subseteq\{1,\dots,m\} \\
\scriptstyle J\ne\emptyset \end{array}}\biggl[p^{|J|p}E_J\max_{{\bf i}_J}
\Bigl(\sum_{{\bf
i}_{J^c}}E_{J^c}h_{\bf i}\Bigr)^p\biggr]\Biggr\}.
\end{eqnarray}
\vskip.1truein
A similar inequality with different constants can be obtained from (2.8'). This
is the most elaborate form we will give to our bounds for
$h\ge 0$ and
$p>1$.

The right hand side of (2.15) for $m=2$ becomes, disregarding constants,
$$E\Bigl(\sum_{i,j}h_{i,j}\Bigr)^p\le
C\max\biggl[
E_1\max_i\Bigl(\sum_jE_2h_{i,j}\Bigr)^p,
~E_2\max_j\Bigl(\sum_iE_1h_{i,j}\Bigr)^p,$$
\vskip-.2truein
$$\mbox{\phantom{aaaaaaaaaaaaaaaaaaaaaaaaaaaaaaaaaaaaaaaaaaa}}
E\max_{i,j}h_{i,j}^p,~t_0^p\biggr].\leqno(2.15')$$
So,  we get the $p$-th moment of the double sum controlled by
moments of partial
maxima of conditional expectations plus a quantile.  The Gin\'e-Zinn (1992)
inequality (for
$m=2$),
$$E\Bigl(\sum_{i,j}h_{i,j}\Bigr)^p\le
C\max\biggl[
E\max_i\Bigl(\sum_jh_{i,j}\Bigr)^p, t_0^p
\biggr],~~p\ge 1,$$
is slightly weaker in appearance than (2.15') (actually, we only published the
result for canonical $U$-statistics, but we applied it
as well to nonnegative variables, for which the proof is the same: see, e.g.,
Gin\'e and Zhang (1996)). For applications of this inequality in the asymptotic
theory of
$U$-statistics see Gin\'e and Zhang (1996), Gin\'e,
Kwapie\'n, Lata\l a and Zinn (1999) and de la Pe\~na and Gin\'e (1999).

\vskip.05truein
\n{\bf Remark.} The constants in the definition of $t_0$ in (2.15)
depend on $p$,
hence, so does $t_0$. This is not the case when $m=1$ (as a consequence of the
improved Hoffmann-J\o rgensen's inequality of Kwapie\'n and Woyczy\'nski -see,
de la Pe\~na and Gin\'e (1999) p. 11-). But in most applications it
does not matter
whether the definition of the quantile depends on $p$.

\vskip.1truein
\n{\it 2.2. Canonical kernels, moments or order $p\ge 2$.} If $\xi_i$
are centered
and independent and $p\ge2$, then, by convexity and the Khinchin-Bonami
inequality
(e.g., de la Pe\~na and Gin\'e, 1999, p. 113), we have
\begin{eqnarray}
2^{-p}E\Bigl(\sum\xi_i^2\Bigr)^{p/2}\!\!\!&\le&\!\!\!
2^{-p}E\Big|\sum\varepsilon_i\xi_i\Big|^p\ \le\  E\Big|\sum\xi_i\Big|^p
\nonumber\\
\!\!\!&\le&\!\!\!
2^pE\Big|\sum\varepsilon_i\xi_i\Big|^p\ \le\ 2^p(p-1)^{p/2}E\Bigl(\sum
\xi_i^2\Bigr)^{p/2},
\end{eqnarray}
where $\varepsilon_i$ are independent identically distributed Rademacher
random variables, independent from $\{\xi_i\}$. Suppose $h_{\bf i}$
is  canonical
for the variables $\{X_i^{(j)}\}$ given in the previous subsection,
that is, suppose
\begin{equation}
E_jh(X_{i_1}^{(1)},\dots,X_{i_m}^{(m)})=0~~{\rm a.s.\ for\ all}~~
j=1,\dots,m,~~ 1\le i_1,\dots,i_m\le n.
\end{equation}
Let $\varepsilon_i^{(j)}$ be an independent Rademacher array independent of
$\{X_i^{(j)}\}$, and set $$\varepsilon_{\bf i}:=\varepsilon_{i_1}^{(1)}\cdots
\varepsilon_{i_m}^{(m)}.$$
Then, recursive application of inequality (2.16) gives
\begin{eqnarray}
2^{-mp}E\Bigl(\sum_{\bf
i}h_{\bf i}^2\Bigr)^{p/2}\!\!\!&\le&\!\!\!2^{-mp}E\Big|\sum_{\bf
i}\varepsilon_{\bf i}h_{\bf i}\Big|^p\ \le\ E\Big|\sum_{\bf i}h_{\bf i}\Big|^p
\nonumber \\
\!\!\!&\le&\!\!\! 2^{mp}E\Big|\sum_{\bf
i}\varepsilon_{\bf i}h_{\bf i}\Big|^p\
 \le\ 2^{mp}(p-1)^{mp/2}E\Bigl(\sum_{\bf
i}h_{\bf i}^2\Bigr)^{p/2}.
\end{eqnarray}
This
inequality reduces estimation of moments of canonical
$U$-statistics to estimation of moments of nonnegative ones (and
conversely), at
least if constants are not an issue. Combined with Proposition 2.1, it
gives the analogue of Rosenthal's inequality for centered variables
and $p>2$, and
if we apply it in conjunction with Corollary 2.2, we obtain the following
inequality:
\vskip.1truein

\n{\bf Proposition 2.4.} {\it If, for $p>2$ and all ${\bf i}\in
\{1,\dots,n\}^m$,
$h_{\bf i}(X_{i_1}^{(1)},\dots,X_{i_m}^{(m)})$ is $p$-integrable and
$E_jh_{\bf i}(X_{i_1}^{(1)},\dots,X_{i_m}^{(m)})=0$ a.s. for all
$j=1,\dots,m$,  then}
\begin{eqnarray}\lefteqn{2^{-mp}\max_{J\subseteq\{1,\dots,m\}}
\biggl[E_J\max_{{\bf
i}_J}
\Bigl(\sum_{{\bf
i}_{J^c}}E_{J^c}h_{\bf i}^2\Bigr)^{p/2}\biggr]
\le E\Big|\sum_{\bf i}h_{\bf i}\Big|^p}\nonumber \\
&&\le
K_m^p\sum_{J\subseteq\{1,\dots,m\}}\biggl[p^{(m+|J|)p/2}E_J\max_{{\bf
i}_J}\Bigl(\sum_{{\bf i}_{J^c}}E_{J^c}h_{\bf
i}^2\Bigr)^{p/2}\biggr]
\end{eqnarray}
{\it for universal constant $K_m<\infty$.}
\vskip.1truein

And, applying Paley-Zygmund with $r=2$, we finally have:

\vskip.1truein
\n{\bf Theorem 2.5.}  {\it Let $h_{\bf i}$  be as in
Proposition 2.4, and let $p>2$. Then, there exist universal
constants $K_m<\infty$
such that, if $t_0$ is defined
as}
$$t_0=\inf\biggl[t\ge 0: \Pr\Bigl\{\Big|\sum_{\bf i}h_{\bf
i}\Big|>t\Bigr\}\le\Bigl(\frac{3}{
4}\Bigr)^{p/(p-2)}\frac{1}{\bigl(2K_m^pp^{mp/2}\bigr)^{1/(p-2)}}\biggr],$$
{\it then}
$$\frac{1}{(4K_mp^{m/2})^{p/(p-2)}}t_0^p\vee
\max\biggl[2^{-mp}\max_{\begin{array}{c}\scriptstyle
J\subseteq\{1,\dots,m\} \\
\scriptstyle J\ne\emptyset \end{array}
}\biggl[E_J\max_{{\bf
i}_J} \Bigl(\sum_{{\bf
i}_{J^c}}E_{J^c}h_{\bf i}^2\Bigr)^{p/2}\biggr]$$
\vskip-.2truein
\begin{equation}\le E\Big|\sum_{\bf i}h_{\bf i}\Big|^p
\mbox{\phantom{aaaaaaaaaaaaaaaaaaaaaaaaaaaaaaaa}}
\end{equation}
\vskip-.2truein
$$\le
2K_m^p\Biggl\{(2p^{m/2})^pt_0^p+\sum_{\begin{array}{c}
\scriptstyle J\subseteq\{1,\dots,m\} \\
\scriptstyle J\ne\emptyset \end{array}}\biggl[p^{(m+|J|)p/2}
E_J\max_{{\bf i}_J}
\Bigl(\sum_{{\bf
i}_{J^c}}E_{J^c}h_{\bf i}^2\Bigr)^{p/2}\biggr]\Biggr\}.$$

\vskip.1truein

If, instead of inequality (2.2), we wish to obtain an analogue
of inequality (2.2'),
that is, if we want to replace the constants at the right hand
side of (2.19) by
$(Kp/\log p)^{mp}$, then we cannot use Khinchin's inequality
and must proceed
directly with an induction as in Proposition 2.1 with the
following change: we must
consider the variables $\sum_{{\bf i}_{J_{m-1}^c}}h_{\bf i}$
as taking values in
$L_2(J_{m-1}^c)$ and apply inequality (1.5) in Kwapie\'n and
Szulga (1991), which
gives Rosenthal's inequality with best constants for centered
independent random
variables in Banach spaces. We skip the details.

\vskip.1truein
\n{\it 2.3. Nonnegative kernels, moments of order $p\le 1$}. It seems
impossible to obtain inequalities as simple as in the previous section for
this case. However, one can still obtain
inequalities that may become useful when combined with Paley-Zygmund.
Here is an analogue of Corollary 2.2 for $h\ge 0$ and $p\le 1$.
The method of proof
is inefficient regarding constants as Hoffmann-J\o rgensen is applied twice at
each step. Hence, constants will not be specified.

\vskip.1truein
\n{\bf Proposition 2.6.} {\it Let  $0< r<p\le 1$, $m<\infty$ and
assume that the
kernels
$h_{\bf i}\ge0$ have integrable $p$-th powers. Then}
\begin{eqnarray}\lefteqn{\max_{J\subseteq\{1,\dots,m\}
}\biggl[E_J\max_{{\bf
i}_J} \Bigl(E_{J^c}\Bigl(\sum_{{\bf
i}_{J^c}}h_{\bf i}\Bigr)^r\Bigl)^{p/r}\biggr]
\le E\Bigl(\sum_{\bf i}h_{\bf i}\Bigr)^p} \nonumber \\
&&\le
K_{r,p,m}\max_{J\subseteq\{1,\dots,m\}}\biggl[E_J\max_{{\bf i}_J}
\Bigl(E_{J^c}\Bigl(\sum_{{\bf
i}_{J^c}}h_{\bf i}\Bigr)^r\Bigl)^{p/r}\biggr],
\end{eqnarray}
{\it where $K_{r,p,m}$ depends only on the parameters $r,\ p,\ m$.}

\vskip.1truein
Note that all the terms in this bound represent a reduction in the
number of sums
{\it except} for the term corresponding to $J=\emptyset$, which
consists of a power
of the
$r$-th moment of a
$U$-statistic of order $m$. We will deal later with this term by means of the
Paley-Zygmund argument.

\vskip.1truein
\n{\bf Proof.} The inequality at the left side of (2.21) follows from H\"older.
Inequality $(H_r)$  is just the right hand side of inequality
(2.21) for
$m=1$ and we can proceed by induction. We still use the notation from
Proposition 2.1. By the induction hypothesis we have
\begin{eqnarray}\label{2ind1} \lefteqn{\mbox{\phantom{aaaaaaa}}
E\Bigl(\sum_{{\bf
i}}h_{{\bf i}}\Bigr)^{p}=E_{m}E_{\{1,\ldots,m-1\}} \Bigl(\sum_{{\bf
i}_{\{1,\ldots,m-1\}}}\sum_{i_{m}}h_{{\bf i}}\Bigr)^{p}}\\
&&\leq K_{r,p,m-1}\! \sum_{J_{m-1}\subset\{1,\ldots,m-1\}}\!
E_{J_{m-1}}E_{m}\max_{{\bf i}_{J_{m-1}}}
\biggl[E_{J_{m-1}^{c}}\Bigl(\sum_{{\bf i}_{J_{m-1}^{c}}}\sum_{i_{m}}
h_{{\bf i}}\Bigr)^{r}\biggr]^{p/r}.\nonumber
\end{eqnarray}
Let us fix $J_{m-1}\subset \{1,\ldots,m-1\}$ and note that, for
fixed $(X_{j}^{(i)})_{i\in J_{m-1}}$, we have
\[ \max_{{\bf i}_{J_{m-1}}}
  E_{J_{m-1}^{c}}\Bigl(\sum_{{\bf i}_{J_{m-1}^{c}}}\sum_{i_{m}}h_{{\bf
i}}\Bigr)^{r}:=
\Big\|\sum_{i_{m}}\tilde{h}_{i_{m}}\Big\|\]
for suitably chosen independent r.v.'s $\tilde{h}_{i_{m}}$
in $l^{\infty}(L^{r})$.
Therefore by ($H_{r}$), which still holds in this space
(as the norm, restricted
to nonnegative vectors, is monotone increasing), we have
\begin{eqnarray} \label{2ind2}\lefteqn{
E_{J_{m-1}}E_{m}\max_{{\bf i}_{J_{m-1}}}
\biggl[E_{J_{m-1}^{c}}\Bigl(\sum_{{\bf i}_{J_{m-1}^{c}}}\sum_{i_{m}}
h_{{\bf i}}\Bigr)^{r}\biggr]^{p/r}
=E_{J_{m-1}}E_{m}\Big\|\sum_{i_{m}}\tilde{h}_{i_{m}}\Big\|^{p/r}}\nonumber\\
&&\leq
  C_{p,r}E_{J_{m-1}}\biggl[E_{m}\max_{i_{m}}\|\tilde{h}_{i_{m}}\|^{p/r}+
\Bigl(E_{m}\big\|\sum_{i_{m}}\tilde{h}_{i_{m}}\big\|\Bigr)^{p/r}\biggr]
\nonumber\\
&&=C_{p,r}\biggl[E_{J_{m-1}\cup\{m\}}\max_{{\bf
i}_{J_{m-1}\cup\{m\}}}\Bigl(E_{J_{m-1}^{c}} \Bigl(\sum_{{\bf
i}_{J_{m-1}^{c}}}h_{{\bf i}}\Bigl)^{r}\Bigr)^{p/r}\nonumber\\
&&\mbox{\phantom{aaaaaaaa}}+E_{J_{m-1}}\Bigl(E_{m}\max_{{\bf
i}_{J_{m-1}}}E_{J_{m-1}^{c}} \Bigl(\sum_{{\bf i}_{J_{m-1}^{c}}}
\sum_{i_{m}}h_{{\bf
i}}\Bigr)^{r}\Bigr)^{p/r}\biggr].
\end{eqnarray}
Now, to estimate the last term, we note that
\begin{eqnarray}\label{2ind3} \lefteqn{E_{J_{m-1}}\biggl[E_{m}\max_{{\bf
i}_{J_{m-1}}}E_{J_{m-1}^{c}} \Bigl(\sum_{{\bf i}_{J_{m-1}^{c}}}
\sum_{i_{m}}h_{{\bf
i}}\Bigr)^{r}\biggr]^{p/r}}\nonumber\\
&&\mbox{\phantom{aaaaaaa}}\leq E_{J_{m-1}}\biggl[E_{J_{m-1}^{c}\cup\{m\}}
  \Bigl(\sum_{{\bf i}}h_{{\bf i}}\Bigr)^{r}\biggr]^{p/r}
\end{eqnarray}
\vskip-.05truein
\n $\leq K_{p/r,1,|J_{m-1}|}$
\vskip-.1truein
$$\times\sum_{J\subset J_{m-1}}E_{J}\max_{{\bf i}_{J}}
\biggl[E_{J_{m-1}\setminus J\cup J_{m-1}^{c}\cup\{m\}}
\Bigl(\sum_{{\bf i}_{J_{m-1}\setminus J\cup J_{m-1}^{c}\cup\{m\}}}
h_{{\bf i}}\Bigr)^{r}\biggr]^{p/r},$$
which follows by the version of Corollary 2.2 for $L^r$ ((2.8'') for $p/r>1$).
Now (\ref{2ind1}), (\ref{2ind2}) and (\ref{2ind3}) complete the induction
step.
\quad$\bbox$

\vskip.05truein
To deal with the term corresponding to $J=\emptyset$ in
Proposition 2.6 we apply
Paley-Zygmund as above, but now with $r<p$ replacing $1<p$.
The conclusion is:

\vskip.1truein

\n{\bf Theorem 2.7}. {\it There is a constant $K_{r,p,m}$ such that
 for  $0< r<p\le1$, $m<\infty$, and
$h_{\bf i}\ge0$ with integrable $p$-th powers, we have}
\begin{eqnarray}\lefteqn{
\frac{1}{(2^{p+1}K_{r,p,m})^{1/(p-r)}}
t_0^p\vee\sum_{\begin{array}{c}\scriptstyle
J\subseteq\{1,\dots,m\} \\
\scriptstyle J\ne\emptyset \end{array}
}\biggl[E_J\max_{{\bf i}_J} \Bigl(E_{J^c}\Bigl(\sum_{{\bf
i}_{J^c}}h_{\bf i}\Bigr)^r\Bigl)^{p/r}\biggr]}\nonumber\\
&&\le E\Bigl(\sum_{\bf i}h_{\bf i}\Bigr)^p \\
&&\le
2K_{r,p,m}\Biggl\{2^{p/r}t_0^p+\sum_{\begin{array}{c}\scriptstyle
J\subseteq\{1,\dots,m\} \\
\scriptstyle J\ne\emptyset \end{array}}\biggl[E_J\max_{{\bf i}_J}
\Bigl(E_{J^c}\Bigl(\sum_{{\bf
i}_{J^c}}h_{\bf i}\Bigr)^r\Bigl)^{p/r}\biggr]\Biggr\},\nonumber
\end{eqnarray}
{\it where}
$$t_0=\inf\biggl[t:\Pr\Bigl\{\sum_{\bf i}h_{\bf i}>t\Bigr\}\le\frac{1}{
2}(2^{p+1}K_{r,p,m})^{-1/(p-r)}\biggr].$$

\vskip.1truein

Hence, the $p$-th moment of a $U$-statistic of order $m$ can be estimated by
partial moments of maxima (or sums) of conditional moments of $U$-statistics
of lower order plus the
$p$-power of a quantile of the original $U$-statistic.

\vskip.1truein
\n{\it 2.4. Canonical kernels and moments of order $1\le p\le 2$,
or kernels $h$
separately symmetric in each of the coordinates and $0<p<1$ }.
The canonical case
reduces to the positive case by means of inequality  (2.18), as before. The
convexity part of inequality (2.18) fails for $p<1$, but in this case,
if $h$ is
symmetric separately in each of the coordinates, we can still
randomize by products
of independent Rademacher variables and recursive application of Khinchin's
inequality still reduces this case to nonnegative $h$. We leave the resulting
statements to the reader in order to avoid repetition.

\vskip.1truein
\n{\it 2.5. Regular (undecoupled) general $U$-statistics.} If $h_{\bf i}({\bf
x})=h_{{\bf i}\circ s}({\bf x}\circ s)$ for any permutation $s$ of
$\{1,\dots,m\}$
and $h_{\bf i}=0$ if ${\bf i}$ has repeated indices, and if the
sequences
$\{X_i^{(j)}:i=1,\dots,n\}$ are independent copies of each other,
then the decoupling
inequalities of de la Pe\~na and Montgomery-Smith (1995), together with the
decoupling inequality for maxima in Hitczenko (1988) in combination
with the previous
inequalities give moment inequalities for the generalized
$U$-statistics
$$\sum_{\bf i}h_{i_1,\dots,i_m}(X_{i_1},\dots,X_{i_m})$$
where $\{X_i\}$ is a sequence of independent random variables,
at the cost of
vastly increasing the numerical constants (see e.g. Gin\'e and Zinn
(1992) for a
similar application of the decoupling inequalities). We omit the resulting
statements.

\vskip.1truein
\n{\it 2.6. Comparison with previous results.} We have already noted, below
the statement of Theorem 2.3, that the inequalities there are better than the
Hoffmann-J\o rgensen
type inequalities for $U$-statistics in Gin\'e and Zinn (1992) in that they
represent a decomposition into simpler quantities. Also, as mentioned in the
Introduction, Ibrahimov and Sharakhmetov (1998, 1999) obtained, except for
constants, Proposition 2.1 and its analogue for canonical kernels for $m=2$
and announced the result for general $m$; the final results in the present
article for
$p>1$ in the nonnegative case (Theorem 2.3) and for
$p>2$ in the canonical case (Theorem 2.5), replacing some sums by maxima and
lower moments by quantiles, seem to be more useful.
As mentioned above, Corollary 2.2 restricted to $m=2$ recovers
inequalities (4.14) in Klass and  Nowicki (1997). The
inequalities in the last mentioned article for nonnegative kernels, $p<1$  and
$m=2$ (the nonconvex case, inequalities (4.13) there) are different from our
inequalities in Theorem 2.7 for $m=2$, although they represent a
similar level of
decomposition of the $p$-th moment of the $U$-statistic.
Basically, the difference is
that they use inverses of truncated conditional moments whereas we use
inverses of
tail probabilites together with partial moments. This can be better
seen by comparing
Hoffmann-J\o rgensen, which is Theorem 2.7 for
$m=1$, with their inequality for $m=1$. The result of Klass and
Nowicki (1997) can
be described  as the iteration of an
inequality that follows from Hoffmann-J\o rgensen, Paley-Zygmund ((2.13)) and
(2.3), as follows. Given $\xi_i$, $i=1,\dots,n$, nonnegative, define
$v_0$ as
\begin{equation}
v_0=\sup\biggl\{v\ge 0: \sum E\Bigl(\frac{\xi_i}{v}\wedge
1\Bigr)\ge 1\biggr\}
\end{equation}
or, what is the same, $v_0$ is the largest number satisfying
\begin{equation}
v_0=\sum E\bigl(\xi_i\wedge v_0\bigr).
\end{equation}
Then, the inequality in question is:

\vskip.1truein
\n{\bf Corollary 2.8.} ({\rm Klass and Nowicki, 1997, Cor. 2.7}) {\it Let
$\xi_i,$ $1=1,\dots,n,$ be independent nonnegative random variables.
Then, for all
$p>0$,}
\begin{equation}
E\Bigl(\sum \xi_i\Bigr)^p\simeq E\max\xi_i^p+v_0^p.
\end{equation}
\vskip.1truein

\n{\bf Proof.} Since
$$\sum E\bigl(\xi_i\wedge\delta_0\bigr)=\sum
E\xi_iI_{\xi_i<\delta_0}+\delta_0\sum\Pr\{\xi_i\ge\delta_0\}\ge\delta_0,$$
it follows that $\delta_0\le v_0$. Therefore, if $p\le 1$,
inequality (2.3) and
the definition of $v_0$ give
$$E\Bigl(\sum\xi_i\Bigr)^p\le\Bigl(\sum E(\xi_i\wedge
v_0)\Bigr)^p+\sum_jE\xi_i^pI_{\xi_i>v_0}\le v_0^p+2E\max_i\xi_i^p.$$
And if $p>1$, Hoffmann-J\o rgensen ($(H)$) and the previous inequality (with
$p=1$) give
$$E\Bigl(\sum\xi_i\Bigr)^p\lessim\Bigl(E\sum\xi_i\Bigr)^p+E\max\xi_i^p\lessim
v_0^p+E\max\xi_i^p.$$
For the reverse inequality, if $p>1$,
$$v_0^p=\Bigl(\sum E(\xi_i\wedge v_0)\Bigr)^p\le E\Bigl(\sum\xi_i\Bigr)^p.$$
And if $p<1$, following the proof of Lemma 2.2 in Klass and
Nowicki (1997), we first
observe that Paley-Zygmund and the first part of this proof give that for some
universal constant
$C$,
\begin{eqnarray*}
\Pr\Bigl\{\sum\xi_i\wedge v_0>\frac{v_0}{2}\Bigr\}\!\!\!&\ge&\!\!\! \frac{1}{
4}\frac{\Bigl(E\sum(\xi_i\wedge v_0)\Bigr)^2}{ E\Bigl(\sum(\xi_i\wedge
v_0)\Bigr)^2}\ =\ \frac{1}{4}\frac{v_0^2}{ E\Bigl(\sum(\xi_i\wedge
v_0)\Bigr)^2} \\
\!\!\!&\ge&\!\!\!\frac{C}{4}\frac{v_0^2}{ E\max (\xi_i\wedge v_0)^2+v_0^2}\
\ge\ \frac{C}{ 8};
\end{eqnarray*}
therefore,
\begin{eqnarray*}
E\Bigl(\sum\xi_i\Bigr)^p\!\!\!&\ge&\!\!\! E\Bigl(\sum(\xi_i\wedge
v_0)\Bigr)^p\\
\!\!\!&\ge&\!\!\! E\biggl[\Bigl(\sum(\xi_i\wedge v_0)\Bigr)^p
I_{\sum(\xi_i\wedge
v_0)>v_0/2}\biggr]\ \ge\ \frac{C}{8}\frac{v_0^p}{2^p}.
\end{eqnarray*}
\n$\bbox$

\vskip.05truein
In fact, if we bound $t_0$ by $t_0^p\le 2
E\Bigl(\sum (\xi_i\wedge t_0)\Bigr)^p$
and apply the above proof to the variables $\xi_i\wedge t_0$,
Hoffmann-J\o rgensen
gives the following seemingly weaker inequality: letting $\tilde{v}_0$  be the
parameter $v_0$ for the smaller variables $\xi_i\wedge t_0$
(note $\tilde{v}_0\le
v_0$), then
$$E\Bigl(\sum \xi_i\Bigr)^p\simeq E\max\xi_i^p+\tilde{v}_0^p.\leqno(2.22')$$

\vskip.2truein

\section {Improved moment inequalities and exponential inequalities for $m=2$}
The right hand side of inequality (2.19) for
$m=1$ is just
\begin{equation}
E\Bigl|\sum\xi_i\Bigr|^p\le
K^p\max\Bigl[p^pE\max\xi_i^p,p^{p/2}\Bigl(\sum
E\xi_i^2\Bigr)^{p/2}\Bigr],\ \ p\ge 2,
\end{equation}
where $\xi_i$ are independent mean zero
random variables. These inequalities were first obtained by Pinelis (1994).
 Part of their interest lie on the fact that they are basically equivalent to
Bernstein's inequality up to constants. Here is how (3.1) (for all $p\ge 2$)
implies Bernstein's inequality up to constants. Assume
$\|\xi_i\|_\infty\le A<\infty$ for all $i$, and set
$C^2=\sum E\xi_i^2$. Then, (3.1) has the form
$$E\Big|\sum\xi_i\Big|^p\le K^p\max\bigl[p^pA^p,~p^{p/2}C^p\bigr],\ \ p\ge 2.$$
 Let
$$p=\frac{x}{KeA}\wedge\Bigl(\frac{x}{KeC}\Bigr)^2$$
for any $x$ for which $p\ge 2$. Then, by Markov's inequality,
(3.1) gives, for these values of $t$,
\[\Pr\Bigl\{\Big|\sum\xi_i\Big|>x\Bigr\}\le\left\{\begin{array}{ll}
\frac{K^pp^pA^p}{x^p}\le e^{-p}       & \mbox{if $p^pA^p\ge p^{p/2}C^p$}\\
&\\
\frac{K^pp^{p/2}C^p}{x^p}\le e^{-p}   & \mbox{otherwise.}
\end{array}
\right.\]
Hence,
\begin{equation}
\Pr\Bigl\{\Big|\sum\xi_i\Big|>x\Bigr\}\le
e^2e^{-p}=e^2\exp\Bigl\{-\frac{x}{KeA}\wedge\Bigl(\frac{x}{
KeC}\Bigr)^2\Bigr\}
\end{equation}
for all $x>0$. Similarly, from the iteration (2.19) of the inequalities (3.1)
we can obtain  exponential inequalities for generalized decoupled
$U$-statistics of any order. However, the inequalities we obtain, while better
than the existing ones, are not of the best kind, as we will see
below. We illustrate this
comment by considering the case
$m=2$. In this case, inequality (2.19) is as follows:
\begin{eqnarray}
E\Bigl|\sum_{i,j}h_{i,j}\Bigr|^p\!\!\!&\le&\!\!\!
K^p\max[p^p\Bigl(\sum_{i,j}Eh_{i,j}^2\Bigr)^{p/2},
p^{3p/2}E_1\max_i\Bigl(\sum_jE_2h_{i,j}^2\Bigr)^{p/2}, \nonumber \\
&&\mbox{\phantom{aaaaaaaa}}p^{3p/2}E_2\max_j
\Bigl(\sum_iE_1h_{i,j}^2\Bigr)^{p/2},
p^{2p}E\max_{i,j}|h_{i,j}|^p\biggr].
\end{eqnarray}
For bounded canonical kernels
$h_{i,j}$ we define
\begin{equation}
A=\max_{i,j}\big\|h_{i,j}\big\|_\infty,\ \
C^2=\sum_{i,j}Eh_{i,j}^2,
\end{equation}
$$B^2=\max\biggl[\Big\|\sum_iE_1h_{i,j}^2(X_i^{(1)},y)\Big\|_\infty,
\Big\|\sum_jE_2h_{i,j}^2(x,X_j^{(2)})\Big\|_\infty\biggr].$$
Then, we can proceed as in
the deduction of (3.2) from (3.1), and easily obtain from (3.3)
that there is a
universal constant
$K$ such that
\begin{equation}
\Pr\biggl\{\Big|\sum_{i,j}h_{i,j}\Big|>x\biggr\}\le
K\exp\biggl\{-\frac{1}{K}\min\Bigl[\frac{x}{
C},\Bigl(\frac{x}{B}\Bigr)^{2/3},\Bigl(\frac{x}{
A}\Bigr)^{1/2}\Bigr]\biggr\}.
\end{equation}
This inequality also holds for regular canonical
$U$-statistics by the decoupling inequalities of de la Pe\~na and
Montgomery-Smith (1995).

Inequality (3.5) is better than the Bernstein type inequality in
Arcones and Gin\'e
(1993) as it is better for
$x\le n^2A$ and the probability is zero for $x\ge n^2A$. Inequality (3.5)
is suboptimal for small values of
$x$, for which the exponent should be a constant times $-x^2$, just as
for chaos variables of order 2 (see Ledoux and Talagrand (1991) and
Lata\l a (1999)).
This suggest that inequality (2.9) is not of the best kind, and can be
improved.

Next we improve the Rosenthal type inequality (2.9) for $m=2$ (that is,
(3.3)) and deduce from it an exponential inequality for canonical
$U$-statistics of order two which does detect the Gaussian portion of the
tail probability.

First we show how Talagrand's (1996) extension of Prohorov's
inequality to empirical processes, actually in Massart's (1999) version,
produces an improved Rosenthal's inequality for empirical processes. Then, we
will use this inequality to estimate the terms resulting from
conditionally applying inequality (3.1) to the $U$-statistic.

To describe Massart's version of Talagrand's inequality we
must establish the setting and define some parameters. Let
$Z_i$ be independent random variables with values in some measurable
space
$(T,{\mathcal T})$, let
$\mathcal F$ be a countable class of measurable real functions on
$T$, and define
$$S:=\sup_{f\in{\mathcal F}}\sum f(Z_i),~~
\sigma^2=\sup_{f\in{\mathcal F}}\sum E(f(Z_i))^2,~~a:=\max_i
\sup_{f\in{\mathcal
F}}\big\|f(Z_i)\bigr\|_\infty.$$
Then,
\begin{equation}\Pr\biggl\{|S|\ge 2E|S|+\sigma\sqrt{8x}+34.5ax\biggr\}\le
e^{-x}
\end{equation}
for all $x>0$.  It follows easily from inequality (3.6) that
\begin{equation}
E|S|^{p}\leq K^{p}\Bigl[(E|S|)^{p}+p^{p/2}\sigma^{p}+
  p^{p}a^{p}\Bigr]
\end{equation}
for some universal constant $K<\infty$ and all $p\ge 1$, in fact, inequality
(3.7) for all $p$ large enough and inequality (3.6) for all $x>0$ are
equivalent up to constants. (We do not plan to keep track of constants in the
derivation below  and, therefore, we refrain from specifying a value for
$K$ in (3.7).)

\vskip.1truein
\n{\bf Proposition 3.1.} {\it Let $\{Z_i\}$ be as above, let $\mathcal F$ be a
countable class of functions such that $Ef^2(Z_i)<\infty$ and $Ef(Z_i)=0$ for
all $i$. Then, in the notation from the previous paragraph,}
\begin{equation}
E|S|^p\le K^p\Bigl[(E|S|)^{p}+p^{p/2}\sigma^{p}+
  p^{p}E\max_i\sup_{f\in{\mathcal F}}\big|f(Z_i)\big|^p\Bigr]
\end{equation}
{\it for all $p\ge 1$, where $K$ is a universal constant.}

\vskip.1truein
\n{\bf Proof.} Set $F:=\sup_{f\in{\mathcal F}}|f|$ and $M^p:=8\cdot
3^pE\max_i|F(Z_i)|^p$. Since the variables $f(Z_i)$ are centered, we can
randomize by independent Rademacher variables $\varepsilon_i$ independent of
the $Z$ variables (at the price of increasing the value of the constant $K$).
Set
$\tilde{S}:=\sup_f\big|\sum\varepsilon_if(Z_i)\big|$. Then,
$$|\tilde{S}|\le\sup_f\big|\sum\varepsilon_if(Z_i)I_{F(Z_i)\le M}\big|+
\sup_f\big|\sum\varepsilon_if(Z_i)I_{F(Z_i)> M}\big|:= S_1+S_2,$$
and notice that, since $ES_1^p\le 2^{p+1}E|S|^p$ (e.g., Lemmas 1.2.6 and
1.4.3 in de la Pe\~na and Gin\'e, 1999), inequality (3.7) gives
$$ES_1^p\le K^p\Bigl[(E|S|)^p+p^{p/2}\sigma^p+p^pM^p\Bigr].$$
To estimate $ES_2^p$ we apply the original Hoffmann-J\o rgensen inequality
(from e.g., Ledoux and Talagrand (1991), (6.9) in page 156) to get
$$ES_2^p\le 2\cdot 3^p\bigl(t_0^p+E\max_iF(Z_i)^p\bigr),$$
where $t_0$ is any number such that $\Pr\{S_2>t_0\}\le(8\cdot 3^p)^{-1}.$
But the choice of $M$ implies that we can take $t_0=0$ because
$$\Pr\bigl\{S_2>0\bigr\}=\Pr\bigl\{\max_iF(Z_i)>M\bigr\}\le{1\over 8\cdot
3^p},$$
proving the proposition.\quad$\bbox$
\vskip.1truein

In what follows we will assume, just as above, that the kernels
$h_{i,j}$, $i,j\leq
n$, are completely degenerate and define
\begin{equation}D=\|(h_{i,j})\|_{L^{2}\rightarrow
L^{2}}:=\sup\biggl\{E\sum_{i,j}h_{i,j}(X_{i}^{(1)},X_{j}^{(2)})
f_{i}(X_i^{(1)})g_{j}(X_j^{(2)})
\end{equation}
$$\phantom{aaaaaaaaaaaaaaaaaaaaaaaaa}:E\sum_{i}f_{i}^{2}(X_i^{(1)})\leq 1,
E\sum_{j}g_{j}^{2}(X_j^{(2)})\leq 1\biggr\}.$$

\vskip.1truein

\n{\bf Theorem 3.2.} {\it There exists a universal constant $K<\infty$  such
that, if $h_{i,j}$ are bounded canonical kernels of two variables for the
independent random variables $X_i^{(1)},X_j^{(2)}$,
$i.j=1,\dots,n$, $n\in\bf{N}$,
 then}
\begin{eqnarray}
E\Big|\sum_{1\le i,j\le
n}\!\!\!\!&&\!\!\!\!h_{i,j}(X_{i}^{(1)},X_{j}^{(2)})\Big|^{p}\le
K^{p}\biggl[p^{p/2}\Bigl(\sum_{i,j} Eh_{i,j}^2\bigr)^{p/2}+p^{p}
\|(h_{i,j})\|_{L^{2}\rightarrow
L^{2}}\nonumber\\
\!\!\!\!&&\!\!\!\!+p^{3p/2}\Bigl[E_1\max_i\Bigl(\sum_jE_2h_{i,j}^2\Bigr)^{p/2}+
E_2\max_j\Bigl(\sum_iE_1h_{i,j}^2\Bigr)^{p/2}\Bigr]\\
&&+p^{2p}E\max_{i,j}|h_{i,j}|^p
\biggr]\nonumber
\end{eqnarray}
{\it for all $p\ge 2$.}

\vskip.1truein

Inequality (3.10) is strictly better than the right hand side inequality in
(2.9) for $m=2$, that is, than (3.3).

\vskip.1truein
\n{\bf Proof.}  Inequality
(3.1) applied conditionally on the variables $X_i^{(1)}$  gives
\begin{equation}
E\Big|\sum_{i,j}h_{i,j}\Big|^p\le K^pE_1\biggl(p^{p/2}
\biggl[\sum_{j}E_2\Bigl(\sum_{i}h_{i,j}\Bigr)^{2}
\biggr]^{p/2}+p^pE_2\sum_j\Big|\sum_ih_{i,j}\Big|^p\biggr).
\end{equation}
To bound the first summand at the right hand side of (3.11) we first notice
that
\begin{eqnarray*}&&\biggl[\sum_{j}E_2\Bigl(\sum_{i}h_{i,j}\Bigr)^{2}
\biggr]^{1/2}\\
&&=\sup\biggl[\sum_{i}E_2\sum_{j}h_{i,j}(X_{i}^{(1)},X_{j}^{(2)})
f_{j}(X_j^{(2)}):
  E\sum_{j}f_{j}^{2}(X_j^{(2)})\leq 1\biggr],
\end{eqnarray*}
where in fact, the sup is taken only over a countable subset of
mean zero vector
functions
$(f_1,\dots,f_n)$ dense in the unit ball of $L_2({\mathcal
L}(X_1^{(2)}))\times\cdots\times L_2({\mathcal L}(X_n^{(2)}))$ for
the seminorm
$|(f_j)_{j\le n}|=\Bigl(\sum Ef_j^2(X_j^{(2)})\Bigr)^{1/2}$. [To see this,
first apply duality in $\ell_2^n$ and then in $L_2(\mathcal{L}(X_j^{(2)}))$
for each $j$.] So we  can apply (3.8) to
$Z_i=(h_{i,j})_{j=1}^n$ with
$f(Z_i)=E_2\sum_{j}h_{i,j}(X_{i}^{(1)},X_{j}^{(2)})f_{j}(X_j^{(2)})$. In this
case, the right hand side terms in (3.8) can be estimated as follows. The
first term:
$$(E|S|)^2\le
E|S|^2=E\biggl[\sum_{j}E_2\Bigl(\sum_{i}h_{i,j}\Bigr)^{2}\biggr]
=E\sum_{i,j}
  h_{i,j}^{2}=C^2.$$
For the second we see that, since, by the previous duality argument,
$$\sum_{i}E_1\Bigl(E_2\sum_{j}h_{i,j}(X_{i}^{(1)},X_{j}^{(2)})
f_{j}(X_j^{(2)})\Bigr)^{2}\leq
  \|(h_{i,j})\|_{L^{2}\rightarrow L^{2}}^2=D^2,$$
it follows that $\sigma\le D$.
The third term:
\begin{eqnarray*}E\max_i\sup_f|f(Z_i)|^p\!\!\!\!&=&\!\!\!\!
E_1\max_i\sup_{E\sum
f_j^2\le1}
\Bigl[E_2
\sum_{j}h_{i,j}(X_{i}^{(1)},X_{j}^{(2)})f_{j}(X_j^{(2)})\Bigr]^p
\\
\!\!\!\!&\le&\!\!\!\! E_1\max_i\sup_{E\sum f_j^2\le1}\Bigl[
  \Bigl(E_2\sum_{j}h_{i,j}^{2}\Bigr)^{1/2}
\Bigl(E\sum_j
f_{j}^{2}\Bigr)^{1/2}\Bigr]^p\\
\!\!\!\!&=&\!\!\!\!
E_1\max_i\Bigl(E_2\sum_{j}h_{i,j}^{2}\Bigr)^{p/2}.
 \end{eqnarray*}
Thus, inequality (3.8) gives
\begin{eqnarray}
p^{p/2}&&\!\!\!\!\!\!\!\!\!\!\!E\Bigl(\sum_{j}E_2\bigl(\sum_{i}h_{i,j}\bigr)^{2}
\Bigr)^{p/2}\\
&&\!\!\!\!\!\!\!\!\!\!\!\le
K^p\Bigl[p^{p/2}C^p+p^pD^p+p^{3p/2}E_1\max_i\Bigl(E_2\sum_{j}h_{i,j}^{2}
\Bigr)^{p/2}\Bigr].\nonumber
\end{eqnarray}
To estimate the second summand at the right hand side
of (3.11), we apply (3.1) once more and obtain
\begin{eqnarray}
&&p^pE_{2}\sum_{j}E_{1}\Big|\sum_{i} h_{i,j}\Big|^{p}\\
&&\phantom{aaaa}\leq K^{p}
  \biggl[p^{3p/2}E_{2}\sum_{j}\Bigl(\sum_{i} E_{1}h_{i,j}^{2}\Bigr)^{p/2}+
   p^{2p}E\sum_{i,j}\big|h_{i,j}\big|^{p}\biggr].\nonumber
\end{eqnarray}
Thus, to complete the proof of the theorem it suffices to replace the sum in
$j$ and the sum in $i,j$ respectively by maxima in $j$ and in $i,j$ on the
terms at the right hand side of this inequality. But this is an easy exercise
of application of inequality (2.6). For completeness sake, here it is.
Applying (2.6) with $\alpha=3$ and $p/2$ instead of $p$, the first term at
the right of (3.13) bounds as:
\begin{eqnarray*}
&&p^{3p/2}E_{2}\sum_{j}\Bigl(\sum_{i}
E_{1}h_{i,j}^{2}\Bigr)^{p/2}\\
&&\phantom{aaaa}\le
2^{1+3p/2}(1+(p/2)^3)\biggl[\Bigl(\frac{p}{2}\Bigr)^{3p/2}
E_2\max_j\Bigl(\sum_{i}
E_{1}h_{i,j}^{2}\Bigr)^{p/2}+C^p\biggr],
\end{eqnarray*}
which produces the conversion of
the sum into a maximum without increasing the order of the multiplicative
constant in front of $C^p$. The second term in (3.13) requires two steps.
First, we apply (2.6) for $p/2$ and $\alpha=4$, conditionally on
$\{X_i^{(1)}\}$:
\begin{eqnarray}
&&~~~~~~~~~p^{2p}E\sum_{i,j}\big|h_{i,j}\big|^{p}\\
&&\le
2^{2p+1}(1+(p/2)^4)E_1\sum_i\biggl[\Bigl(\frac{p}{2}\Bigr)^{2p}
E_2\max_j|h_{i,j}|^p+\Bigl(\sum_jE_2h_{i,j}^2\Bigr)^{p/2}\biggr].\nonumber
\end{eqnarray}
We apply (2.6) with respect to $E_1$, for $p/2$ and  $\alpha=0$, to the second
term at the right hand side of (3.14) and we obtain the bound
$$2^{2p+3}(1+(p/2)^4)\biggl[E_1\max_i
\Bigl(\sum_jE_2h_{i,j}^2\Bigr)^{p/2}+C^p\biggr],$$
which is in terms of some of the quantities appearing at the right hand side
of (3.10) and with coefficients of lower order. As for the first term at the
right of (3.14), we apply (2.6) with respect to
$E_1$, again for
$p/2$ and
$\alpha=4$, and get it bounded by
$$2^{4p+2}(1+(p/2)^4)^2\biggl[\Bigl(\frac{p}{2}\bigr)^{2p}
E\max_{i,j}|h_{i,j}|^p+E_2\Bigl(\sum_iE_1\max_jh_{i,j}^2\Bigr)^{p/2}\biggr].$$
Here the first term coincides with the last one in (3.10), and the second is
dominated by
$$K^pE_2\Bigl[\sum_j\Bigl(\sum_iE_1h_{i,j}^2\Bigr)\Bigr]^{p/2}.$$
Applying inequality $(R_1)$ with respect to $E_2$ this is in turn dominated by
$$K^p\Bigl(\frac{p}{2}\Bigr)^{p/2}E_2\sum_j
\Bigl(\sum_iE_1h_{i,j}^2\Bigr)^{p/2} +K^pC^p,$$
and the first summand has alredy been handled above
(first term at the right of (3.13)). Collecting terms we obtain inequality
(3.10).\quad$\bbox$

\vskip.1truein

Theorem 3.2 gives the following moment inequality and
exponential bound for bounded kernels.

\vskip.1truein
\n{\bf Theorem 3.3.} {\it There exist universal constants $K<\infty$ and
$L<\infty$ such that, if $h_{i,j}$ are bounded canonical kernels of two
variables for the independent random variables $X_i^{(1)},X_j^{(2)}$,
$i.j=1,\dots,n$, and if
$A$, $B$, $C$, $D$ are as defined in (3.4) and (3.9), then}
\begin{equation}
E\Big|\sum_{1\le i,j\le n}h_{i,j}(X_{i}^{(1)},X_{j}^{(2)})\Big|^{p}\leq
K^{p}\Bigl[p^{p/2}C^{p}+p^{p}D^{p}+
  p^{3p/2}B^{p}+p^{2p}A^{p}\Bigr]
\end{equation}
{\it for all $p\ge 2$ and, equivalently,}
\begin{eqnarray}\Pr\biggl\{\bigg|\sum_{i,j\leq
n}h_{i,j}(X_{i}^{(1)},X_{j}^{(2)})\bigg|\!\!\!&\geq&\!\!\! x\biggr\}
\nonumber\\
\!\!\!&\leq&\!\!\! L\exp\biggl[-\frac{1}{ L}\min\Bigl(\frac{x^{2}}{
C^{2}},
   \frac{x}{D},\frac{x^{2/3}}{B^{2/3}},
   \frac{x^{1/2}}{A^{1/2}}\Bigr)\biggr]
\end{eqnarray}
{\it for all $x>0$.}

\vskip.1truein
The moment inequality is immediate from Theorem 3.2 and the equivalence with
the exponential inequality follows just like (3.2) follows from (3.1) in one
direction, and, in the other, by integration of tail probabilities.

Next we comment on the exponential inequality.
For comparison purposes, let
$h_{i,j}(X_i^{(1)},X_j^{(2)})=g_ig_j'x_{i,j}$ with $g_i,g'_j$
independent standard
normal. In this case,
$$C^2=\sum_{i,j} x_{i,j}^2~~~{\rm and}~~~D=\sup
\Bigl\{\sum_{i,j}u_iv_jx_{i,j}:~\sum
u_i^2\le1,\sum v_j^2\le1\Bigr\}$$
and the Gaussian chaos inequality in  Lata\l a (1999)
yields the existence of  universal constants $0<k<K<\infty$ such that
$$\Pr\Bigl\{\Big|\sum_{i,j}h_{i,j}\Big|\ge K(Cx^{1/2}+Dx)\Bigr\}\le e^{-x}$$
and
$$\Pr\Bigl\{\Big|\sum_{i,j}h_{i,j}\Big|\ge k(Cx^{1/2}+Dx)\Bigr\}
\ge k\wedge e^{-x}.$$
By the central limit theorem for canonical $U$-statistics, this
implies that the
coefficients of
$x^2$ and $x$ in (3.16) are correct (except for $K$). It is
natural to have terms in smaller powers of $x$ in (3.16) e.g., by
comparison with Bernstein's inequality for sums of independent random
variables. In
fact, the term in $x^{1/2}$ cannot be avoided, at least up to
logarithmic factors.
To see this, consider the product $V$ of two independent
centered Poisson variables
with parameter 1, which is the limit in law of $V_n=\sum_{i,j\le
n}X_i^{(n)}Y_j^{(n)}$ where
$X_i^{(n)}$ and $Y_j^{(n)}$ are centered Bernoulli random variables
with parameter
$p=1/n$; then, for large $x$, the tail probabilities of $V$ are of
the order of
$\exp{(-x^{1/2}\log x)}$, and therefore, so are those of $V_n$ for large $n$.
Also, note that the term in $x^{2/3}$ in the exponent corresponds, up to
logarithmic factors, to the tail probabilities of the product of
two independent
random variables, one normal and the other centered Poisson.

If $X,Y,X_i^{(1)},X_j^{(2)}$ are i.i.d., $h_{i,j}=h$ for all $i,j$ and $h$ is
completely degenerate, then the parameters defined by (3.4) and (3.8) become:
$$A=\|h\|_\infty,~~B^2=
n\bigl(\|E_Yh^2(x,Y)\|_\infty+\|E_Xh^2(X,y)\|_\infty\bigr),~~C^2=n^2Eh^2$$ and
\begin{eqnarray*}D\!\!\!&=&\!\!\!n\sup\Bigl\{Eh(X,Y)f(X)g(Y):~Ef^2(X)\le
1,Eg^2(Y)\le 1\Bigr\}\\
\!\!\!&:=&\!\!\!n\|h\|_{L_2\mapsto L_2},
\end{eqnarray*}
where $\|h\|_{L_2\mapsto L_2}$ is the norm of the operator of $L_2({\mathcal
L}(X))$ with kernel $h$. Then, inequalities (3.15) and (3.16) become:
\vskip.1truein

\n{\bf Corollary 3.4.} {\it Under the above assumptions,
there exist universal
constants $K<\infty$, $L<\infty$ such that, for all $n\in\bf N$ and $p\ge 2$,}
\begin{eqnarray}
&&E\Big|\sum_{i,j\le n}h(X_i^{(1)},X_j^{(2)})\Big|^p\le
K^p\Bigl[p^{p/2}n^p(Eh^2)^{p/2}+p^pn^p\|h\|_{L_2\mapsto
L_2}^p\nonumber\\
&&\phantom{aaaaaaaaaaaaaa}+p^{3p/2}n^{p/2}
\bigl(\|E_Yh^2\|_\infty+\|E_Xh^2\|_\infty\bigr)^{p/2}
+p^{2p}\|h\|_\infty^p\Bigr]
\end{eqnarray}
{\it and}
\begin{eqnarray}&&\Pr\biggl\{\Big|\sum_{i,j\le n}h(X_i^{(1)},X_j^{(2)})\Big|\ge
x\biggr\}
\le K\exp\biggl[-\frac{1}{K}\min\biggl(\frac{x^2}{
n^2Eh^2},\nonumber\\
&&\mbox{\phantom{aaaaaaaa}}\frac{x}{n\|h\|_{L_2\mapsto
L_2}},\frac{x^{2/3}}
{\bigl[n(\|E_Yh^2\|_\infty+\|E_Xh^2\|_\infty)\bigr]^{1/3}},
\frac{x^{1/2}}{\|h\|_\infty^{1/2}}\biggr)\biggr].
\end{eqnarray}

Inequality (3.18) provides an analogue of Bernstein's
inequality for degenerate
$U$-statistics of order 2: note that inequalities (3.15), (3.16), (3.17) and
(3.18) can all be `undecoupled' using the result of de la Pe\~na and
Montgomery-Smith's (1995). It should also be noted that this exponential
inequality for canonical $U$-statistics  is strong enough to imply the
sufficiency part of the
law of the iterated logarithm for these objects: this
can be seen by applying it to
the kernels $h_n$ in Steps 7 and 8 of the proof of Theorem 3.1
in Gin\'e, Kwapie\'n, Lata\l
a and Zinn (1999) (and using some of the computations there for the
parameters $C$ to $D$). Neither inequality (3.5) nor any of the previously
published inequalities for $U$-statistics can do this.

\vskip.15truein
\n{\bf Acknowledgement.} We thank Stanislaw Kwapie\'n for several useful
conversations.

\vskip.15truein
\centerline{\bf References}
\vskip.2truein \baselineskip=10pt
\parskip=4pt

{\smc\n Arcones, M. and Gin\'e, E.} (1993). Limit theorems for $U$-processes.
{\it Ann. Probab.} {\bf 21} 1494-1542.\par
{\smc\n de la Pe\~na, V. and Gin\'e, E.} (1999). {\it Decoupling:
>From Dependence to
Independence}. Springer-Verlag, New York.
\par
{\smc\n de la Pe\~na, V. and Montgomery--Smith, S.} (1995).
Decoupling inequalities for the tail probabilities of multivariate
$U$-statistics.
{\it Ann. Probab.} {\bf 23} 806-816.\par
{\smc\n Figiel, T.; Hitczenko, P.; Johnson, W.B.; Schechtman, G.; and Zinn, J.}
(1997). Extremal properties of Rademacher functions with applications to the
Khintchine and Rosenthal inequalities. {\it Trans. Amer. Math. Soc.} {\bf 349}
997-1027.
\par
{\smc\n Gin\'e, E.; Kwapie\'n, S.; Lata\l a, R.; and Zinn, J.} (1999). The LIL
for canonical $U$-statistics of order two. To appear.\par
{\smc\n Gin\'e, E. and Zhang, C.-H.} (1996). On integrability in the LIL for
degenerate $U$-statistics. {\it J. Theoret. Probab} {\bf 9} 385-412.\par
{\smc\n Gin\'e, E. and Zinn, J.} (1983). Central limit theorems and weak
laws of
large numbers in certain Banach spaces. {\it Zeits. Wahrsch. v. Geb.} {\bf 62}
323-354.\par
{\smc\n Gin\'e, E. and Zinn, J.} (1992). On Hoffmann-J\o rgensen's
inequality for
$U$-processes. {\it Probability in Banach Spaces 8} 80-91. Birkh\"auser,
Boston.\par
{\smc\n Hitczenko, P.} (1988). Comparison of moments for tangent sequences
of random variables. {\it Probab. Th. Rel. Fields} {\bf 78} 223-230.\par
{\smc\n Johnson, W. B.; Schechtman, G.; and Zinn, J.} (1985). Best constants
in moment inequalities for linear combinations of independent and exchangeable
random variables. {\it Ann. Probab.} {13} 234-253.\par
{\smc\n Ibragimov, R. and Sharakhmetov, Sh.} (1998). Exact bounds on the
moments of symmetric statistics. In: Abstracts of the 7-th Vilnius Conference
on Probability Theory and Mathematical Statistics/ 22nd European Meeting of
Statisticians, pp. 243-244. Vilnius.\par
{\smc\n Ibragimov, R. and Sharakhmetov, Sh.} (1999). Analogues of Khintchine,
Marcinkiewicz-Zygmund and Rosenthal inequalities for symmetric statistics.
{\it Scand. J. Statist.} {\bf 26} 621-623.\par
{\smc\n Kahane, J.-P.} (1968). {\it Some Random Series of Functions}.
Heath, Lexington, Massachusetts.\par
{\smc\n Klass, M. and Nowicki, K.} (1997). Order of magnitude bounds for
expectations of $\Delta_2$ functions of nonnegative random bilinear forms and
generalized $U$-statistics. {\it Ann. Probab.} {\bf 25} 1471-1501.\par
{\smc\n Kwapie\'n, S. and Szulga, J.} (1991). Hypercontraction
methods in moment
inequalities for series of independent random variables in normed spaces.
{\it Ann. Probab.} {\bf 19} 369-379.\par
{\smc\n Kwapie\'n, S. and Woyczy\'nski, W.} (1992). {\it Random Series
and Stochastic
Integrals: Single and Multiple}. Birkh\"auser, Boston.\par
{\smc\n Lata\l a, R.} (1997). Estimation of moments of sums of
independent random
variables. {\it Ann. Probab.} {\bf 25} 1502-1513.\par
{\smc\n Lata\l a, R.} (1999). Tails and moment estimates for some type
of chaos. {\it Studia Math.} {\bf 135} 39-53.\par
{\smc\n Lata\l a, R. and Zinn, J.} (1999).
Necessary and sufficient conditions for the strong law of large numbers for
$U$-statistics. {\it Ann. Probab.}, to appear.\par
{\smc\n Ledoux, M. and Talagrand, M.} (1991). {\it Probability in Banach
Spaces:
Isoperimetry and Processes}. Springer, New York.\par
{\smc\n Massart, P.} (1999). About the constants in Talagrand's concentration
inequalities for empirical processes. {\it Ann. Probab.}, to appear.\par
{\smc\n Pinelis, I.} (1994). Optimum bounds for the distributions
of martingales
in  Banach spaces. {\it Ann. Probab.} {\bf 22} 1679-1706.\par
{\smc\n Talagrand, M.} (1996). New
concentration inequalities in product spaces. {\it Invent. Math.} {\bf 126}
505-563.\par
{\smc\n Utev, S. A.} (1985). Extremal problems in moment inequalities. In: {\it
Limit Theorems in Probability Theory}, Trudy Inst. Math.,
Novosibirsk, 56-75 (in
Russian).\par

\vskip.2truein

\begin{tabbing}
Department of Mathematics and Statisticslalala\=Institute of Mathematics\kill
Evarist Gin$\acute{\mbox{e}}$\>Rafa{\l} Lata\l a\\
Department of Mathematics\>Institute of Mathematics\\
and Department of Statistics\>Warsaw University\\
University of Connecticut\>Banacha 2\\
Storrs, CT 06269\>02-097 Warszawa\\
USA\>Poland\\
gine@uconnvm.uconn.edu\>rlatala@mimuw.edu.pl
\end{tabbing}

\vskip.2truein

\begin{tabbing}
aaaaaaaaaaaaaaaaaaaaa\=Department of
mathematical Sciences\kill
\>Joel Zinn\\
\>Department of Mathematics\\
\>Texas A\&M University\\
\>College Station, TX 77843\\
\>jzinn@math.tamu.edu
\end{tabbing}

\end{document}